\documentclass{amsart}
\usepackage{amssymb}
\usepackage{bm}
\usepackage{amsmath}
\usepackage{mathrsfs}
\usepackage[dvips]{graphicx}
\usepackage{graphics}
\usepackage{subfigure}
\usepackage{flafter}

\newtheorem{theorem}{Theorem}[section]

\newtheorem{definition}[theorem]{Definition}

\newtheorem{lemma}[theorem]{Lemma}

\newtheorem{remark}[theorem]{Remark}

\begin{document}

\title[Lattice on Caratheodory Extension]{Construction of A Lattice on the completion space of an algebra and an isomorphism to its Caratheodory Extension}

\author{Jun Tanaka}
\address{University of California, Riverside, USA}
\email{juntanaka@math.ucr.edu, yonigeninnin@gmail.com, junextension@hotmail.com}

\author{Peter F. Mcloughlin}
\address{Crafton Hills College, USA}
\email{pmcloughlin@aol.com}

\keywords{Caratheodory Extension Theorem, Lattice Theory, Measure Theory}
\subjclass[2000]{Primary: 28A12, 28B20}
\date{August, 10, 2008}

\begin{abstract}
In this paper, we will show how one is able to construct a lattice on the completion of an algebra and to obtain an isomorphism to its Caratheodory Extension. In addition, it will be shown that the lattice form a $ \sigma$-algebra and a complete Heyting algebra of countable type.
\end{abstract}

\maketitle

\section{\textbf{Introduction}}
The Caratheodory Extension Theorem is a crucial part of any advanced
Real Analysis course; the process extends an algebra to a $
\sigma$-algebra and a measure on an algebra to a measure on a $
\sigma $-algebra. This is a powerful tool in both measure theory and
statistics (see for example, \cite{fleming},\cite{stochastic}). A Boolean Algebra structure on the Caratheodory Extension was discussed in several papers (e.g. Kolmogorov, \cite{Boolean1}, Coquand and Palmgren \cite{Boolean2}). Furthermore, the Caratheodory Extension is complete with a pseudometric $ d ( . , .) = \mu^{\ast} (. \triangle . ) $ where $\mu^{\ast} ( . ) $ is the outer measure on the power set P(X) (Please refer to \cite{Dun}). The Caratheodory Extension has various rich structures which are open to further investigation.

In this paper, we will show how the Caratheodory Extension process
is intimately related to the metric completion process. In
particular, it will be shown how one is able to construct a lattice on the completion of an algebra and to obtain an isomorphism
to its Caratheodory Extension. The author believes this method to be new, intuitive, and constructive. Especially noteworthy is that, in Definition \ref{d:9} of this paper, we give a precise definition of $\vee^{\infty}$ on the completion space $ ( \bar{d}, \overline{\Omega}$) of an algebra $\Omega$ and show that the lattice is a complete Heyting algebra of countable type.

Let X be an arbitrary nonempty set, and let $ \mu$ be a finite measure on an algebra $ \Omega \subset \mathbf{P}(X)$.
Note that $ d ( . , .) = \mu^{\ast} (. \triangle . ) $ defines a pseudometric on
$ \Omega$. Denote by $(\bar{d} , \overline{\Omega})$ the completion
of $ (d, \Omega)$. Also, let $S$ be the set of all Cauchy sequences
in $(d, {\Omega})$. By the metric completion procedures, we know $
\bar{d} ( \{ \overline{B^{\alpha_{1}}_n} \}, \{
\overline{B^{\alpha_{2}}_n} \}) = \lim d  ( {B^{\alpha_{1}}_n} ,
{B^{\alpha_{2}}_n}) $ where $ \{  {B^{\alpha_{i}}_n} \} $ is in $S$
with $ \{ \overline{B^{\alpha_{i}}_n} \} = \{ \{B_n \}  \in S : \lim
d (B_n, B_n ^{\alpha_{i}} ) = 0 \} $ for $ i  = 1 $or 2.

\begin{definition}\label{d:1}
Let $ \bar{\mu} ( \{ \overline{B^{\alpha}_{n}} \})  = \bar{d} ( \{
\overline{B^{\alpha}_{n}} \} , \{ \bar{\phi} \} ) = \lim \mu (
{B^{\alpha}_{n}}  )$.
\end{definition}

\begin{definition}\label{d:2}
If $ B _n \in \Omega$, then $ \{ B_n \} $ is called a $ \mu$-Cauchy
sequence if $\mu (B_n \triangle B_m ) \rightarrow 0 $ as n,m $\rightarrow \infty$.

\end{definition}

\begin{definition}\label{d:3}
For A, B in $\mathbf{P}(X) $, A = B a.e. if $\mu^{\ast} ( A  \triangle B )  =0. $

\end{definition}

In Section \ref{Se:Completion}, we will show that a lattice structure can be naturally defined on $
\overline{\Omega} $ which makes it a $ \sigma$-algebra. Let $ \widetilde{\mathbf{S}}$ = $\{S \in \mathbf{P}(X) \ |$ $ \exists \ \mu$-Cauchy sequence  $\{B_{n}   \}   $ s.t. $ \lim \mu^{\ast} ( B_{n}   \triangle S)  =0   \}$. In \cite{Junext}, we proved that $ \widetilde{\mathbf{S}}$ is a $\sigma$-algebra where, for any $\mu$-Cauchy sequence $ \{ B_{n} \}$ such that $ \lim \mu^{\ast} ( B_{n}   \triangle S)  =0 $, the measure $\widetilde{\mu} (S) $ on $ \widetilde{\mathbf{S}}$ is defined as $\widetilde{\mu} (S) $ = $\lim \mu ( B_{n}   )   $. In addition, we proved that $  \widetilde{\mu}$ is a countably additive measure on $ \widetilde{\mathbf{S}}$. Thus, ($\widetilde{\mu}$, $ \widetilde{\mathbf{S}}$) is a measure space. We showed that the Caratheodory Extension of $\Omega$ can be expressed as the set of limit points of $\mu$-Cauchy sequences under the pseudometric $d(A,B)$ = $\mu^{\ast}(A \triangle B)$. Moreover, when the measure is a $\sigma$-finite measure, we obtained an equivalent expression of the Caratheodory Extension, $ \widetilde{\mathbf{S}}$ = $\{S \in \mathbf{P}(X) $ $ | $  $\exists \ \mu$-Cauchy sequence  $\{B_{n}   \}   $ s.t. $ \lim \mu^{\ast} ( B_{n}   \triangle S)  =0   \}$. Theorem 2 in \cite{Junext} shows that $E$ is a measurable set iff  $E$ is in $ \widetilde{\mathbf{S}}$. Thus, the measure space ($\widetilde{\mu}$, $ \widetilde{\mathbf{S}}$) agrees with the Caratheodory Extension when $\mu$ is a finite measure. Moreover, it shows that measurable sets are exactly limit points of $\mu$-Cauchy sequences. The $\sigma$-finite case follows from the finite case. In Section \ref{Se:Iso}, we will define a $\sigma$-algebra lattice isomorphism between two $\sigma$-algebras, and define a map $F :$ $\overline{\Omega}   \rightarrow \mathbf{P}(X) $ given by $F ( \overline{\{B_{n}\}} )  = B $ where $ \lim \mu^{\ast} ( B_{n}   \triangle B  )  =0   $. We will show that $F$ is an isometry and a $\sigma$-algebra lattice isomorphism between the completion $\overline{\Omega}$ and the Caratheodory Extension of $\Omega$ under the equivalence relation $\sim$ defined as A $\sim$ B iff $\mu^{\ast}(A \triangle B)$ =0.

\section{\textbf{Main Definitions and Notations}}

Our notation agrees with that in \cite{Junext}. For emphasis, let X be a nonempty set. Throughout the paper, unless otherwise stated, $ \mu$ will be a finite measure on an algebra $ \Omega \subset \mathbf{P}(X)$. Unless otherwise stated, $\{ B^\alpha_n\}$ will be a $ \mu$-Cauchy
sequence.

In addition, for each $ \mu$-Cauchy
sequence $\{ B_n\}$, a capital letter corresponds to the limit point of the $ \mu$-Cauchy
sequence. For instance,
$ \lim \mu^{\ast} ( A_{n}   \triangle A  )  =0   $, $ \lim \mu^{\ast} ( B_{n}   \triangle B  )  =0   $, $ \lim \mu^{\ast} ( Y_{L}   \triangle Y  )  =0   $. There exist such limit points in $\mathbf{P}(X)$ since the Caratheodory Extension is complete under the pseudometric d $($\cite{Dun}$)$.

We recall that by a $\sigma$-algebra of subsets of X, we mean a family $\Omega$ of subsets such that (1) every countable union of sets in $\Omega$ is in $\Omega$, and (2) $A^C$ is in $\Omega$ whenever A is.

\begin{definition}\label{d:5}
Define a map $F :$ $\overline{\Omega}   \rightarrow \mathbf{P}(X) $ given by $F ( \overline{\{B_{n}\}} )  = B $

where $ \lim \mu^{\ast} ( B_{n}   \triangle B  )  =0   $. Note that such a B always exists in $\mathbf{P}(X)$ since the Caratheodory Extension is complete $($\cite{Dun}$)$.

\end{definition}

\begin{definition}\label{d:6}{\textbf{$\sigma$-algebra lattice isomorphism}}

Suppose lattices $($X, $\vee_X$,$\wedge_X$ $)$ and $($Y, $\vee_Y$,$\wedge_Y$ $)$ are $\sigma$-algebras, and H: X $\rightarrow  $ Y is a one-to-one, onto well-defined map. Then H is called a $\sigma$-algebra lattice isomorphism if
\[
\begin{aligned}
H(\cdot  \vee_X \cdot ) = & H(\cdot ) \vee_Y  H(\cdot )  , \  \  \ \  \ \  \  H(\vee^{\infty}_{X } \ \cdot ) & = &  \vee^{\infty}_{Y}   H(\cdot )  ,  \\
H(\cdot  \wedge_X \cdot ) = & H(\cdot ) \wedge_Y H(\cdot )  ,\ \ \  \ \ \ \ H(\cdot^{\textbf{C}} ) & = & H(\cdot )^{\textbf{C}} .
\end{aligned}
\]

\end{definition}

We will show that F is a $\sigma$-algebra lattice isomorphism between $\overline{\Omega}$ and $ \widetilde{\mathbf{S}}_{\diagup_{\sim} }$ where $\sim$ is an equivalence relation in the sense of Definition \ref{d:3}.

\section{\textbf{Preliminaries}}

In this section, we shall briefly review the well-known facts about lattice theory (e.g. Birkhoff \cite{Birk}, Iwamura \cite{Iwamura}), propose an extension lattice, and investigate its properties, as well as some lemmas from \cite{Junext}. A nonempty set L is called a lattice if L is closed under the operations of meet ($\wedge  $) and join ($\vee  $). It is denoted by (L,$\wedge  $,$\vee  $) or simply L. If it satisfies, in addition, the distributive law, then it is called a distributive lattice. For two lattices L and $L'$, a bijection from L to $L'$, which preserves lattice operations is called a lattice isomorphism, or simply an isomorphism. If there is an isomorphism from L to $L'$, then L is called lattice-isomorphic with $L'$, and we write L $\cong$ $L'$. We write x $\leq $ y if x $\wedge$ y = x or, equivalently, if x $\vee$ y = y. L is called complete if, for any subset A of L, L contains the supremum $\vee $A and the infimum $\wedge$A, with respect to the above order. A complete lattice L includes the maximum and minimum elements, which are denoted by I and O, or 1 and 0, respectively. A distributive lattice is called a Boolean algebra or a Boolean lattice, if, for any element x in L, there exists a unique complement $x^C \in L$ such that x $\vee$ $x^C$ = 1 and $x \wedge x^C$ = 0.
Let L be a lattice and $\cdot^{c}$: L $\rightarrow $ L be an operator. Then $\cdot^{c}$ is called a lattice complement in L if the following conditions are satisfied:

(1) $\forall x$ $\in$ L, $x^C \wedge x$ = 0, $x^C \vee x$ = I,

(2) $\forall$ x, y $\in$ L, x $\leq$ y $\Rightarrow$ $x^C \geq y^C$,

(3) $\forall$ x $\in$ L, $(x^C)^C$ = x.

\begin{definition}

A complete lattice L is called a lattice $\sigma$-algebra or simply a $\sigma$-algebra if the following conditions are satisfied:

(1) $\forall$ A $\in$ L, $A^C \in $L where $ \cdot^C$ is a lattice complement,

(2) $\forall$ $A_i \in$ L, $\vee^{\infty} A_i \in$ L.

\end{definition}

\begin{definition}\label{cha:1}

A complete lattice is called a complete Heyting algebra $($cHa$)$, if
\[
\begin{aligned}
\vee_{i \in I} \ ( x_{i} \wedge y ) =  (\vee_{i \in I} \  x_{i} ) \wedge y
\end{aligned}
\]
holds for $\forall x_{i} , y \in L$ $(i \in I)$; where I is an index set of arbitrary cardinal number. In the case that the cardinality of I is countable, it is called a complete Heyting Algebra of countable type.
\end{definition}

It is well-known that for a set E, $| P(E) |$= $2^{|E|}$. The set of all subsets of E is a Boolean algebra.

\begin{lemma}\label{le:1}
If $ \{B_n\} $ is a $ \mu$-Cauchy sequence, then $ \{ \mu (B_n) \} $
is a Cauchy sequence of real numbers.
\end{lemma}

\begin{proof}
Note $ B_n \subseteq (B_n \triangle B_m) \cup B_m $ and $ B_m
\subseteq (B_n \triangle B_m) \cup B_n $ implies $ \mu(B_n) \leq \mu
(B_n \triangle B_m) + \mu (B_m) $ and $ \mu (B_m ) \leq \mu (B_n
\triangle B_m) + \mu (B_n) $. Hence, $| \mu (B_n) - \mu (B_m ) |\leq
\mu (B_n \triangle B_m )$.
\end{proof}

\begin{lemma}\label{le:2}
$ d(A _1 \cup A_2, A_3 \cup A_4 ) \leq d (A_1, A_3) + d (A_2 , A_4 )
$ for any $ A _i \in \Omega $.
\end{lemma}

\begin{proof}
The proof is in \cite{Junext}
\end{proof}

\begin{lemma}\label{le:3}
$\{B^\alpha _n \cup B^\gamma _n \} $ and $ \{ ( B ^\alpha _n )^C \}
$ are $ \mu$-Cauchy sequences.
\end{lemma}

\begin{proof} By Lemma $\ref{le:2}$,
\[d(B^\alpha _ n \cup B^\gamma _n , B ^\alpha _m \cup B^\gamma _m ) \leq d (B^\alpha _n, B^\alpha _m ) + d (B ^\gamma _n, B^\gamma _m ).
\]
\[
\begin{aligned}
& \mbox{Moreover}, \, \,  d ((B^\alpha _n )^C , (B^\alpha _m )^C ) = \mu ( (B^\alpha _n )^C \triangle (B^\alpha _m )^C ) =  \mu (B^\alpha _n  \Delta B^\alpha _m ) = d (B^\alpha _n ,
B^\alpha _m ).\\  & \mbox{Hence, the claim follows}.
\end{aligned}
\]
\end{proof}

\section{\textbf{The Completion Space $ ( \bar{d}, \overline{\Omega}$) and the Lattice (L,$\wedge  $,$\vee  $)} }\label{Se:Completion}

In this section, we introduce the completion of ($d$,$\Omega$) where $\Omega$ is an algebra of subsets of X.

\begin{definition}\label{d:7}
Let $E_\alpha = \{ \overline{ B^\alpha_n}\} $ and $E_A = \overline{ \{
A \} } $ when $ A \in \Omega $ and $\{  B^\alpha_n \} $ is a sequence in $\Omega$.
\end{definition}

\begin{definition}\label{d:8}
Define $ E_\alpha \vee E_\gamma = \{ \overline{B^\alpha _n \cup B
^\gamma _n } \} $, $ E_\alpha \wedge E_\gamma = \{ \overline{B^\alpha _n \cap B
^\gamma _n } \} $, and $ E^C _\alpha  = \{ \overline{(B^\alpha _n )
^C }\} $.
\end{definition}

\begin{remark}\label{n:1}
$ E _\alpha \wedge E_\gamma = E _\phi $ iff $ \bar{\mu} (E_\alpha
\wedge E_\gamma) = 0 $ iff $ \lim \mu (B^\alpha _n \cap B ^\gamma
_n ) = 0 $. Moreover, $ E_\alpha \subset E_\gamma $ iff $ \bar{\mu}
( E _\alpha \wedge E ^C_\gamma ) = 0 $. The containment may be written in the usual lattice sense; namely, $ E_\alpha \wedge E_\gamma = E_\alpha$.
\end{remark}

\begin{lemma}\label{le:4}
The map $ \vee : \overline{\Omega} \times \overline{\Omega} \rightarrow
\overline{\Omega} $ defined by $ \vee (E_\alpha \times E_\gamma ) =
E_\alpha \vee E_\gamma = \{ \overline{B^\alpha _n \cup B^\gamma
_n} \}   $ is well-defined.
\end{lemma}

\begin{proof} Suppose $\{ \overline{B^\alpha _n} \} = \{ \overline{A^\alpha _n} \} $ and $ \{ \overline{A^\gamma _n} \} = \{ \overline{B^\gamma _n} \}$.

Then $ \lim d (B^\alpha _n, A ^\alpha _n) = \lim d (B^\gamma _n,
A^\gamma _n)=0 $. Also, by Lemma 3, $ \{B^\alpha _n \cup
B^\gamma _n \} $ and $ \{ A^\alpha _n \cup A^\gamma _n\} $ are $
\mu$-Cauchy sequences.

Hence, by Lemma $\ref{le:2}$, $ d (B^\alpha _n \cup B^\gamma _n, A^\alpha
_n\cup A^\gamma _n) \leq d(B^\alpha _n , A^\alpha _n) + d( B^\gamma
_n, A^\gamma _n) $ and the claim follows.
\end{proof}

\begin{lemma}\label{le:5} The map $ \cdot {}^C : \overline{\Omega} \rightarrow \overline{\Omega}  $ defined by \[ (\{ \overline{B^\alpha _n} \}) ^C = \overline{ \{ ( B^\alpha _n)^C \} } \]
is well-defined.
\end{lemma}

\begin{proof}
Suppose $ \{ \overline{B^\alpha _n} \} = \{ \overline{A^\alpha _n}
\}$. Then $ \lim d (B ^\alpha _n , A ^\alpha _n ) = 0$. Now, $d ((
B^\alpha _n)^C, (A^\alpha _n)^C )= \mu (( B^\alpha _n )^C \triangle
(A^\alpha _n)^C )= \mu (B^\alpha _n \triangle A^\alpha _n) = d (
B^\alpha _n , A^\alpha _n ) $ which implies $ \{ ( B^\alpha _n )^C
\} \sim \{(A^\alpha _n)^C\}$.
\end{proof}

\begin{lemma}\label{le:4}
The map $ \wedge : \overline{\Omega} \times \overline{\Omega} \rightarrow
\overline{\Omega} $ defined by $ \wedge (E_\alpha \times E_\gamma ) =
E_\alpha \wedge E_\gamma = \{ \overline{B^\alpha _n \cap B^\gamma
_n} \}   $ is well-defined.
\end{lemma}

\begin{proof}
The proof follows immediately from the two previous lemmas.
\end{proof}

\begin{lemma}
$ \cdot {}^C $ in Lemma \ref{le:5} is a lattice complement.
\end{lemma}

\begin{proof}
The proof is obvious.
\end{proof}

\begin{lemma}\label{le:6} $ \bar{\mu} (E_\alpha \vee E_\gamma ) \leq \bar{\mu} (E _{\alpha} ) + \bar{\mu} (E_\gamma) $ for $ E_\alpha , E _\gamma \in \overline{\Omega} $.
\end{lemma}

\begin{proof}
$E_\alpha = \{ \overline{B^\alpha _n} \} $ and $ E_\gamma = \{
\overline{B^\gamma _n} \} $ implies  $E_\alpha \vee E_\gamma = \{
\overline{B^\alpha _n \cup B^\gamma _n} \} $. Hence $ \bar{\mu} ( E
_\alpha \vee E _\gamma )= \lim \mu ( B^\alpha _n \cup B^\gamma _n
)  \leq \lim [ \mu (B^\alpha _n) + \mu (B^\gamma _n)] = \lim \mu
(B^\alpha _n) + \lim \mu ( B^\gamma _n) \leq \bar{\mu}(E_\alpha ) +
\bar{\mu} (E _\gamma)$.
\end{proof}

\begin{lemma}\label{le:7}
If $ E _{\alpha_{1}} $ and $ E _{\alpha_{2}}$ are disjoint, then $
\bar{\mu} (E _{\alpha_{1}} \vee E _{\alpha_{2}} ) = \bar{\mu} (E
_{\alpha_{1}} ) + \bar{\mu} (E _{\alpha_{2}} )$.
\end{lemma}

\begin{proof}
 By Remark $\ref{n:1}$, we must have $ \lim \mu ( B ^{\alpha_{1}} _n \cap B ^{\alpha_{2}} _n)=0$. Note, $B ^{\alpha_{i}} _n = (B ^{\alpha_{i}} _n \cap B ^{\alpha_{j}} _n) \ \cup \ ( B ^{\alpha_{i}} _n \cap (B ^{\alpha_{j}} _n) ^C )$ implies $ \mu (B ^{\alpha_{i}} _n) = \mu ( B ^{\alpha_{i}} _n  \cap B ^{\alpha_{j}} _n )+ \mu (B ^{\alpha_{i}} _n \cap (B ^{\alpha_{j}} _n)^C )$ and $\lim \mu (B ^{\alpha_{i}} _n ) = \lim \mu (B ^{\alpha_{i}} _n \cap (B ^{\alpha_{j}} _n )^C )$ for $ \{i , j \} = \{1, 2\} $.

Moreover,  $ \mu ( B ^{\alpha_{1}} _n   \cup B ^{\alpha_{2}} _n  ) =
\mu ( \ ( B ^{\alpha_{1}} _n \cap B ^{\alpha_{2}} _n ) \ \cup \ (B ^{\alpha_{1}} _n  \cap ( B ^{\alpha_{2}} _n  )^C ) \ \cup \ (B ^{\alpha_{2}} _n   \cap (B ^{\alpha_{1}} _n
 )^C ) \ ) = \mu (B ^{\alpha_{1}} _n
 \cap B ^{\alpha_{2}} _n
)+ \mu (B ^{\alpha_{1}} _n
 \cap (B ^{\alpha_{2}} _n
 )^C )+ \mu (B ^{\alpha_{2}} _n  \cap (B ^{\alpha_{1}} _n
)^C) $ implies $ \bar{\mu} (E _{\alpha_{1}} \vee E _{\alpha_{2}})
= \lim \mu
 (B ^{\alpha_{1}}_n  \cup B ^{\alpha_{2}}_n)= \lim \mu (B ^{\alpha_{1}}_n)+ \lim \mu (B ^{\alpha_{2}}_n)= \bar{\mu} ( E_{\alpha_{1}} ) + \bar{\mu} (E_{\alpha_{2}}) $.
\end{proof}

\begin{theorem}\label{le:8}
$ \overline{\Omega} $ is an algebra and $ \bar{\mu} $ is a measure
on $ \overline{\Omega} $.
\end{theorem}

\begin{proof} By Lemmas $\ref{le:3}$, $\ref{le:4}$  and $\ref{le:5}$,  $ \overline{\Omega}$ is an algebra. Moreover,  the fact that $\bar{\mu}$ takes on values in $\mathbb{R} ^+ $ along with Lemmas 6 and 7 imply $ \bar{\mu} $ is a measure on  $\overline{\Omega} $.
\end{proof}

\begin{remark}\label{n:2}  $ \bar{d} ( E_{\alpha},  E _\gamma ) = \lim d (B ^{\alpha}_n  ,   B ^{\gamma}_n ) = \lim [ \mu (B ^{\alpha}_n  \cap (B ^{\gamma}_n )^C )+ \mu (B ^{\gamma}_n \cap (B ^{\alpha}_n )^C)] = \lim \mu (B ^{\alpha}_n  \cap (B ^{\gamma}_n )^C)+\lim \mu (B ^{\gamma}_n  \cap (B ^{\alpha}_n )^C) = \bar{\mu} ( E_{\alpha} \wedge (E_{\gamma})^C)+ \bar{\mu} ( E_{\gamma} \wedge (E_{\alpha}) ^C) = \bar{\mu} (E_{\alpha} \triangle E_{\gamma}) $.
\end{remark}

\begin{lemma}\label{le:1}

Let $E_{i}$ =  $\overline{\{B^{i}_{n}\}}$  $\in \overline{\Omega}$ for i $\geq$ 1 and by following the proof of Lemma 8 in \cite{Junext}, construct $Y_{L} $ =  $\cup_{i=1}^{N_{L}} B^{i}_{K_{L}}$ for each $L$ such that

\[
\mu^{\ast}(\cup_{i=1}^{\infty}   S_{i}  \triangle  \cup_{i=1}^{N_{L}} B^{i}_{K_{L}}     ) \leq  \mu^{\ast}(\cup_{i=N_{L}+1}^{\infty}   S_{i}) + \mu^{\ast}(\cup_{i=1}^{N_{L}}   S_{i}  \triangle  \cup_{i=1}^{N_{L}} B^{i}_{K_{L}}     ) < \frac{1}{L}.
\]

Then
E := $\overline{\{ Y_{L}   \}}  $ satisfies the following conditions;

\[
\begin{aligned}
1. & \vee ^n _{i=1} E _{i}  \leq E \ \ \mbox{for all} \  n, \\
2. &\lim \bar{\mu} (E \wedge ( \vee ^n_{i=1} E _{i}
)^C ) = 0   .
\end{aligned}
\]

In particular, E is uniquely determined in $\overline{\Omega}$.

\end{lemma}

\begin{proof}
Note that $E_{i}$ =  $\overline{\{B^{i}_{n}\}}$ = $\overline{\{B^{i}_{K_{L}}\}}$.

$( \vee_{i=1}^{n}   E_{i}  ) \wedge \overline{\{ Y_{L}   \}} $  =  $ \overline{\{  \cup_{i=1}^{n} B^{i}_{K_{L}}   \cap  Y_{L}         \}}   $ = $ \overline{\{  \cup_{i=1}^{n} B^{i}_{K_{L}}     \}}   $ = $\vee_{i=1}^{n}   E_{i}   $ for any n.

Let $N_{L} > n$.
\[
\begin{aligned}
\mu( \cup_{i=1}^{N_{L}} B^{i}_{K_{L}}  \cap  (\cup_{i=1}^{n} B^{i}_{K_{L}})^{\textbf{C}}   ) = & \mu( \cup_{i=1}^{N_{L}} B^{i}_{K_{L}}  \triangle  \cup_{i=1}^{n} B^{i}_{K_{L}}     )  = \mu^{\ast}( \cup_{i=1}^{N_{L}} B^{i}_{K_{L}}  \triangle  \cup_{i=1}^{n} B^{i}_{K_{L}}     ) \\  \leq & \mu^{\ast}(\cup_{i=1}^{\infty}   S_{i}  \triangle  \cup_{i=1}^{N_{L}} B^{i}_{K_{L}}     ) + \mu^{\ast}(\cup_{i=1}^{\infty}   S_{i}  \triangle  \cup_{i=1}^{n} B^{i}_{K_{L}}     ).
\end{aligned}
\]

This implies that $\lim  \overline{\mu}( \overline{\{ Y_{L}   \}}  \wedge ( \vee_{i=1}^{n}     E_{i}) ^{\textbf{C}}  )   $ = 0.

Suppose $E, E^\prime \in \overline{\Omega} $ both satisfy condition 1 and 2 above. Then $  \vee ^n _{i=1} E _{i}  \leq E \wedge E^\prime \leq E$ for all n. Hence, we must have $ E = E \wedge E ^\prime $. Similarly, we must have $ E^\prime = E \wedge E^\prime $. Therefore E is unique.

\end{proof}

\begin{lemma}\label{le:14}
Let $E_{i}$ =  $\overline{\{B^{i}_{n}\}}$  $\in \overline{\Omega}$ for i $\geq$ 1 and $ E_{i} \wedge E_{j}  = E_\phi $ for $ i
\neq j$. Then $ \bar{\mu} (E) = \sum ^\infty _{i=1} \bar{\mu}
(E_{i}) $ where E is defined as in Lemma \ref{le:1}.

\end{lemma}

\begin{proof}
Suppose that $E_{i}$ =  $\overline{\{B^{i}_{n}\}}$  $\in \overline{\Omega}$ for i $\geq$ 1 and $ E_{i} \wedge E_{j}  = E_\phi $ for $ i \neq j$. Then it follows that $ \bar{\mu} (\vee ^n _{i=1} E _{i}) = \sum^k_{i=1} \bar{\mu} ( E_{i} ) \leq \bar{\mu} (E) $ for all $k$.

Hence, we must have
\begin{equation}\label{Eq:1}
\sum^\infty_{i=1} \bar{\mu} ( E_{i } ) \leq
\bar{\mu} (E).
\end{equation}

Now since  $ \bar{\mu} (E) = \bar{\mu}
(\vee ^n _{i=1} E _{i}) + \bar{\mu} (E \wedge (\vee ^n _{i=1} E _{i} )^C )$  for all n and $\lim \bar{\mu} (E \wedge ( \vee ^n_{i=1} E _{i}
)^C ) = 0  $ ,  we must have $ \bar{\mu}
(E) = \sum^\infty _{i=1} \bar{\mu} ( E_{i} ) $ by (\ref{Eq:1}).

\end{proof}

\begin{definition} \label{d:9} For $ {E}_{i}$ in $ \overline{\Omega} $, we define $ \vee ^\infty _{i=1} {E}_{i} = E $ where $E$ is as in Lemma $\ref{le:1}$.
\end{definition}

\begin{theorem}\label{le:10}
$ ( \bar{\mu} , \overline{\Omega} )$ is a measure space.
\end{theorem}

\begin{proof}
This follows by Theorem $\ref{le:8}$, $\ref{le:1}$, and $\ref{le:14}$, and Definition $\ref{d:9}$.
\end{proof}

\section{\textbf{The Isomorphism} }\label{Se:Iso}

In this section, we introduce an isomorphism $F :$ $\overline{\Omega}   \rightarrow \mathbf{P}(X) $ given by $F ( \overline{\{B_{n}\}} )  = B $.

\begin{lemma}

$F$ : $\overline{\Omega}   \rightarrow  \widetilde{\mathbf{S}}_{\diagup_{\sim} } $ is a well-defined map.
\end{lemma}

\begin{proof}
Suppose that $\overline{\{A_{n}\}} =  \overline{\{B_{n}\}}$.

There exist A and B in $\mathbf{P}(X)$ such that $ \lim \mu^{\ast} ( A_{n}   \triangle A  )  =0   $  and  $ \lim \mu^{\ast} ( B_{n}   \triangle B  )  =0   $.

$ \mu^{\ast} (  A  \triangle B  )  \leq  $ $ \mu^{\ast} (  A  \triangle  A_{n}  )  +   \mu^{\ast} (   A_{n}  \triangle B_{n}  ) +   \mu^{\ast} (  B_{n}  \triangle B )  $ by the triangle inequality.

By taking limits on both sides, $ \mu^{\ast} (  A  \triangle B  ) = 0$.

Thus, A = B a.e.. Therefore, F is well-defined.

\end{proof}

\begin{theorem}\label{t:3}
F is an isometry between $\overline{\Omega}$ and $ \widetilde{\mathbf{S}}_{\diagup_{\sim} }$.
\end{theorem}

\begin{proof}
First, we show F is onto $ \widetilde{\mathbf{S}}$. Let $X \in \widetilde{\mathbf{S}}$. Then there exists a $\mu$-Cauchy sequence  $\{B_{n}   \}   $ such that $ \lim \mu^{\ast} ( B_{n}   \triangle X)  =0   $. Thus $F ( \overline{\{B_{n}\}} )    = X $ a.e.. Therefore, F is onto.

Second, we will show F preserves the metric. Let $\overline{\{A_{n}\}}$, $ \overline{\{B_{n}\}}$ $\in $  $\overline{\Omega}$. Letting A and B be as before, we have
\[
  \lim \mu (   A_{n}  \triangle B_{n}  )  =    \mu^{\ast} (  A  \triangle B  )   .
\]
Therefore, $\overline{d}( \overline{\{A_{n}\}}, \overline{\{B_{n}\}} ) $ = $\lim \mu (   A_{n}  \triangle B_{n}  ) $ = $\mu^{\ast} (  A  \triangle B  )     $

=  $\mu^{\ast} (  F ( \overline{\{A_{n}\}} )   \triangle   F ( \overline{\{B_{n}\}} )  )     $ = $d(  F ( \overline{\{A_{n}\}} )  ,   F ( \overline{\{B_{n}\}} )     ) $. Thus, F preserves the metric.

Lastly, we will show that F is one-to-one. Let $  F ( \overline{\{A_{n}\}} )  ,   F ( \overline{\{B_{n}\}} )    \in \widetilde{\mathbf{S}}$ such that $  F ( \overline{\{A_{n}\}} )  =  F ( \overline{\{B_{n}\}} ) $ a.e..

Then $ A  =  B   $ a.e. implies $\mu^{\ast} (  A  \triangle B  ) = 0$. Then, as in the proof of F being onto,  $\lim \mu (   A_{n}  \triangle B_{n}  ) $ = $\mu^{\ast} (  A  \triangle B )     $. Thus $\overline{\{A_{n}\}}$ = $\overline{\{B_{n}\}}$ and F is one-to-one. Therefore, F is an isometry between $\overline{\Omega}$ and $ \widetilde{\mathbf{S}}_{\diagup_{\sim} }$.

\end{proof}

\begin{theorem}\label{t:4}
F is a $\sigma$-algebra lattice isomorphism between $\overline{\Omega}$ and $ \widetilde{\mathbf{S}}_{\diagup_{\sim} }$.
\end{theorem}

\begin{proof}
We already showed that F is a one-to-one, onto map in Theorem $\ref{t:3}$.

$F ( \overline{\{A_{n}\}}  \vee  \overline{\{B_{n}\}} )$ = $F ( \overline{\{A_{n} \cup B_{n}\}} )$ = $F ( \overline{\{A_{n}\}} )  \cup F(  \overline{\{B_{n}\}} )$ a.e. since $\lim \mu (   A_{n} \cup B_{n} \ \ \triangle \ \ A \cup B  ) $ = 0.

Let $\overline{\{B_{n}\}} \in \overline{\Omega}$.

Then,
\[
F ( \overline{\{B_{n}\}}^{\textbf{C}} )  =    F ( \overline{\{  ( B_{n})^{\textbf{C}}   \}} )   =   B^{\textbf{C}} \ \ a.e..
\]
Thus, $F(\cdot^{\textbf{C}} )$ =  $F(\cdot )^{\textbf{C}} $ in $ \widetilde{\mathbf{S}}_{\diagup_{\sim} }$.

Similarly, $F(\cdot  \wedge \cdot )$ = $F(\cdot^{\textbf{C}} \vee \cdot^{\textbf{C}} )^{\textbf{C}} $ =  $[ F(\cdot^{\textbf{C}} ) \cup  F(\cdot^{\textbf{C}} ) ]^{\textbf{C}} $= $F(\cdot ) \cap F(\cdot )  $.

Let $E_{\alpha_{i}}$ $\in \overline{\Omega}$ for i $\geq$ 1 and $E_{\alpha_{i}}$ = $\overline{\{B^{\alpha_{i}}_{n}\}}$.

Then for each i, there exists an $S_{i}   \in   \widetilde{\mathbf{S}}  $ such that $ \lim \mu^{\ast} (B^{\alpha_{i}}_{n} \triangle S_{i})  =0  $.

Now suppose we have $\{ Y_{L} \}$ in the same manner as in Lemma $\ref{le:1}$. By design, $\{ Y_{L} \}$ converges to $\cup_{i=1}^{\infty}   S_{i} $. Then $\overline{ \{ Y_{L} \}  }$ = $\vee_{i=1}^{\infty}   E_{\alpha_{i}}  $ by Lemma $\ref{le:1}$. Now we have

\[
F( \vee_{i=1}^{\infty}   E_{\alpha_{i}}  )  = F(\vee_{i=1}^{\infty}   \overline{\{B^{\alpha_{i}}_{n}\}}  )  = F ( \overline{ \{ Y_{L} \}  } ) = Y  .
\]
Since $\lim \mu^{\ast} ( Y_{L} \triangle Y  )  =0  $ and $\lim \mu^{\ast} ( Y_{L} \triangle \cup^{\infty}_{i=1} S_{i}  )  =0  $, we have $Y   $ = $\cup^{\infty}_{i=1} S_{i}  $ a.e.. In addition, $\cup^{\infty}_{i=1} S_{i}  $ = $\cup^{\infty}_{i=1} F(  \overline{\{B^{\alpha_{i}}_{n}\}}  )  $.

Thus,
\[
F( \vee_{i=1}^{\infty}   E_{\alpha_{i}}  )  = \cup^{\infty}_{i=1} F( E_{\alpha_{i}} )  .
\]
Therefore, the claim follows.

\end{proof}

\begin{theorem}\label{t:countable}
$( \overline{\Omega}$,$\wedge  $,$\vee  )$ is a complete Heyting algebra of countable type. More precisely, it is the case that the cardinality of I in Definition $\ref{cha:1}$ is countable.
\end{theorem}

\begin{proof}
Since  $ \widetilde{\mathbf{S}}$ is a complete Heyting algebra of countable type, the proof follows from Theorem $\ref{t:4}$.
\end{proof}

\section{\textbf{Conclusion}}
Theorems $\ref{t:3}$ and  $\ref{t:4}$ show that the completion of $\Omega$ is isometric and $\sigma$-algebra lattice isomorphic to $ \widetilde{\mathbf{S}}_{\diagup_{\sim} }$. Thus, by the conclusion in \cite{Junext} the completion of $\Omega$ is isometric to, as well as $\sigma$-algebra isomorphic to, the Caratheodory Extension under the equivalence relation $\sim$ . An important note is that the isomorphism shows that Caratheodory Extension Theorem and the Metric Completion Process are essentially the same, differing only on measure zero sets. In addition, the completion $\overline{\Omega}$ of $\Omega $ is a complete Heyting algebra of countable type by Theorem \ref{t:countable}. The $\sigma$-finite case follows from the finite case.

\section{\textbf{Acknowledgement}}
The first author would like to thank Professors who gave him very professional advice and suggestions, which he truly believes improved the presentation of this paper.

\nocite{*}
\bibliographystyle{aip}

\bibliography{uctest}

\end{document}